# FINITE SAMPLE APPROXIMATION RESULTS FOR PRINCIPAL COMPONENT ANALYSIS: A MATRIX PERTURBATION APPROACH[1]

### By Boaz Nadler

#### *Weizmann Institute of Science*

Principal component analysis (PCA) is a standard tool for dimensional reduction of a set of $n$ observations (samples), each with $p$ variables. In this paper, using a matrix perturbation approach, we study the nonasymptotic relation between the eigenvalues and eigenvectors of PCA computed on a finite sample of size $n$, and those of the limiting population PCA as $n \to \infty$. As in machine learning, we present a finite sample theorem which holds with high probability for the closeness between the leading eigenvalue and eigenvector of sample PCA and population PCA under a spiked covariance model. In addition, we also consider the relation between finite sample PCA and the asymptotic results in the joint limit $p, n \to \infty$, with $p/n = c$. We present a matrix perturbation view of the "phase transition phenomenon," and a simple linear-algebra based derivation of the eigenvalue and eigenvector overlap in this asymptotic limit. Moreover, our analysis also applies for finite $p, n$ where we show that although there is no sharp phase transition as in the infinite case, either as a function of noise level or as a function of sample size $n$, the eigenvector of sample PCA may exhibit a sharp "loss of tracking," suddenly losing its relation to the (true) eigenvector of the population PCA matrix. This occurs due to a *crossover* between the eigenvalue due to the signal and the largest eigenvalue due to noise, whose eigenvector points in a random direction.

**1. Introduction.** Principal component analysis (PCA) is a standard tool for dimensionality reduction, applied in regression, classification and many other data analysis tasks in a variety of scientific fields [17, 20]. PCA finds orthonormal directions with maximal variance of the data and allows its

Received February 2007; revised February 2008.

[1]Supported by the Lord Sieff of Brimpton memorial fund and by the Hana and Julius Rosen fund.

low-dimensional representation by linear projections onto these directions. This dimensionality reduction is a typical preprocessing step in many classification and regression problems.

Assuming the given data is a finite and random sample from a (generally unknown) distribution, an interesting theoretical and practical question is the relation between the sample PCA results computed from finite data and those of the underlying population model. In this paper we consider a *spiked covariance model* for which the underlying data is low-dimensional but each sample is corrupted by additive Gaussian noise, and analyze the following question: how close are the directions and eigenvalues computed by sample PCA to the limiting (and generally unknown) directions and eigenvalues of the population model and how do these quantities depend on the number of samples $n$, the dimensionality $p$ and the noise level $\sigma$.

Many works have studied the convergence of eigenvalues and eigenvectors of sample PCA to those of population PCA under various settings. In general, the different results regarding convergence can be divided into two main types: (i) asymptotic results of classical statistics, where $p$ is fixed and $n \to \infty$, starting with the classical works of Girshick [12], Lawley [23] and Anderson [1], who assumed multivariate Gaussian distributions, up to more recent work which relaxed some of these assumptions; see [2, 17] and references therein; (ii) modern "large $p$, large $n$" statistical results, where the joint limit $p, n \to \infty$ is considered while the ratio $p/n = c$ is kept fixed. In the statistical physics literature we mention Hoyle and Rattray [14], Reimann, Van den Broeck and Bex [32], Watkin and Nadal [37], Biehl and Mietzner [5] and references therein, whereas in the statistics community see the recent works by Johnstone [18], Baik and Silverstein [4], Debashis [31], Onatski [29] and many references therein for older works. However, it seems that the most relevant case, that of *explicit* approximation bounds between the eigenvalues and eigenvectors computed in a *finite* setting ($p, n$ finite) and those of the infinite setting ($p$ fixed, $n = \infty$), as well as estimates of the distributions of these quantities (again for finite and fixed $p, n$), are not covered by these approaches.

In the present work we address this problem using a combination of matrix perturbation theory and concentration of measure bounds on the norm of noisy Wishart matrices. This paper contains two main contributions. First we present *probabilistic approximation results* regarding the difference between the leading eigenvalue and eigenvector of sample PCA and population PCA for *fixed $p$ and $n$*, under a spiked covariance model with a single component. The second contribution of this paper is a matrix perturbation view of the phase transition for the leading eigenvalue and eigenvector, when both $p, n$ are large. We present a simple linear-algebra based proof of the asymptotics of the leading eigenvalue and eigenvector in the limit $p, n \to \infty$. Second, we present an interesting connection between this asymptotic value and a classical result by Lawley. This observation leads to two additional propositions, one regarding the limiting eigenvalues for a more general spiked covariance model, and one regarding the spectral norm of a noisy Wishart matrix with nonidentity diagonal covariance matrix. These results may be useful for inference on the number of significant components under more general settings of heteroscedastic correlated noise. Finally, for finite $p, n$ we show that while there is no deterministic phase transition at a specific fixed value of $p/n$, either as a function of noise level or as a function of sample size, the eigenvector of sample PCA may exhibit a sharp "loss of tracking," suddenly losing its relation to the (true) eigenvector of the population PCA matrix. This occurs due to a *crossover* between the eigenvalue due to the signal and the largest eigenvalue due to noise, whose eigenvector points in a random direction.

Matrix perturbation theory, including eigenvalue and eigenvector perturbation bounds, as well as the structure of eigenvalues and eigenvectors of arrowhead matrices, play a key role in the analysis of both finite sample PCA and the asymptotic limit $p, n \to \infty$. Perturbation theory and concentration of measure results on the norm of noisy Wishart matrices are key ingredients in the analysis of the finite sample case, and also provide novel insight into the origins of the phase transition in the joint limit $p, n \to \infty$. In the statistics literature, matrix perturbation theory is typically used to bound remainder terms in asymptotic limits. In the context of PCA, for example, in [10] Eaton and Tyler used a perturbation bound by Wielandt to present a simple derivation of asymptotic results as $n \to \infty$ for eigenvalues of random symmetric matrices, but did not consider the nonasymptotic case. In [35], Stewart introduced a general framework of stochastic perturbation theory to analyze the effects of random perturbations on the eigenvalues of finite matrices, whereas perturbation bounds for the singular value decomposition were considered in [36]. Within the context of the spiked covariance model, matrix perturbation theory was used in [4, 19, 31], though these works considered mainly the asymptotic limit $p, n \to \infty$.

From a practical point of view, our results show that when $p, n$ are large, and specifically when $p \gg n$, the true signal directions may be drowned by noise since for finite $p \gg n$, eigenvector reconstruction errors behave as $\sigma \sqrt{p/n}$, as also predicted by the asymptotic analysis of Johnstone and Lu [19]. A similar phenomenon also occurs in linear discriminant analysis [6], and for various multivariate regression algorithms such as partial least squares and classical least squares [26]. All these works emphasize the importance of regularization, feature selection or low-dimensional representation prior to learning, and hint that global methods may not be optimal for dimensional reduction or as a preprocessing step prior to regression and classification of high-dimensional noisy data, specifically if there is a priori knowledge regarding its sparsity or smoothness. The results and approach presented in this paper can also be used to develop methods to determine the number of components in a linear mixture (spiked covariance) model [22].

The paper is organized as follows. In Section 2 we present the spiked covariance model and our main results. The results for finite $p, n$ are proven in Sections 3 and 4. A matrix perturbation view of the phase transition phenomenon in the joint limit as $p, n \to \infty$, as well an analysis of the spiked covariance model under more general models of noise and some finite sample examples appear in Section 5.

## 2. Model, assumptions and main results.

2.1. *Notation.* Univariate random variables are denoted by lowercase letters, as in $u$, their realizations are denoted by $u^\nu$ and their expectation by $\mathbb{E}\{u\}$. Column vectors are denoted by boldface lowercase letters, as in $\mathbf{x}$, their transpose is $\mathbf{x}^T$, their $j$th component is $x_j$ the dot product between two vectors is $\langle \mathbf{x}, \mathbf{z} \rangle$, and the Euclidean ($L_2$) norm of $\mathbf{x}$ is $\|\mathbf{x}\|$. Matrices are denoted by uppercase letters, as in $A$. The identity matrix of order $p$ is $I_p$, and the spectral norm of a matrix $A$ is $\|A\|$.

2.2. *The spiked covariance model.* We consider a spiked covariance model where each data sample $\mathbf{x}$ has the form

$$(2.1) \qquad \mathbf{x} = \sum_{i=1}^{k} u_i \mathbf{v}_i + \sigma \xi,$$

where $\{u_i\}_{i=1}^k$ are random variables, typically called components or latent variables, $\{\mathbf{v}_i\}_{i=1}^k \subset \mathbb{R}^p$ are the corresponding fixed (and typically unknown) response vectors, $\xi$ is a multivariate Gaussian noise vector with identity covariance matrix and $\sigma$ is the level of noise. Equation (2.1) is an error-in-variables linear mixture model, commonly assumed in many different problems and applications, including independent component analysis (ICA) [15], signal processing, and in the analysis of spectroscopic data, where it is known as Beer's law [25, 27].

We denote by $\Sigma$ the $p \times p$ population covariance matrix corresponding to the observations $\mathbf{x}$,

$$(2.2) \qquad \Sigma = \mathbb{E}\{\mathbf{x}\mathbf{x}^T\},$$

and by $S_n$ the sample covariance matrix corresponding to the $n$ i.i.d. observations $\{\mathbf{x}^\nu\}_{\nu=1}^n$,

$$(2.3) \qquad S_n = \frac{1}{n} \sum_{\nu=1}^{n} \mathbf{x}^\nu (\mathbf{x}^\nu)^T.$$

Assuming that all $k$ vectors $\mathbf{v}_j$ are orthogonal and that all $k$ random variables in (2.1) are uncorrelated with zero mean, unit variance and finite fourth moment, it follows that the largest $k$ eigenvalue/eigenvector pairs of $\mathbf{\Sigma}$ are given by $\{(\|\mathbf{v}_j\|^2 + \sigma^2, \mathbf{v}_j)\}_{j=1}^k$. PCA approximates the eigenvalues and eigenvectors of the unknown $\Sigma$ by those of $S_n$. In particular, the top eigenvectors of $S_n$ correspond to orthogonal directions of maximal variance of the *observed* data. The question at hand, then, is how close are the largest eigenvalues and corresponding eigenvectors of $S_n$ to those of $\Sigma$.

In this paper we present a matrix perturbation view of this problem. For simplicity we consider the case of a single component $(k = 1)$. A similar though more complicated analysis can be carried out for the case of $k$ components. We thus consider $n$ samples from the model

$$(2.4) \qquad \mathbf{x} = u\mathbf{v} + \sigma\xi,$$

where we assume that the random variable $u$ has zero mean, unit variance and finite fourth moment. Without loss of generality, we further assume that the vector $\mathbf{v}$ is in the direction of $\mathbf{e}_1 = (1, 0, \ldots, 0)$, for example, $\mathbf{v} = \|\mathbf{v}\|\mathbf{e}_1$. Finally, since $u$ has zero mean, we neglect in our analysis the initial mean centering step typically done in PCA.

The population covariance matrix corresponding to this one-component model is given by

$$(2.5) \qquad \Sigma = \begin{pmatrix} \|\mathbf{v}\|^2 & 0 & \cdots & 0 \\ 0 & 0 & \cdots & 0 \\ \vdots & & \ddots & \vdots \\ 0 & 0 & \cdots & 0 \end{pmatrix} + \sigma^2 I_p.$$

Its largest eigenvalue is $\|\mathbf{v}\|^2 + \sigma^2$ with corresponding eigenvector $\mathbf{e}_1$, and all other eigenvalues equal $\sigma^2$.

2.3. *Results for finite $p, n$.* To study the relationship between the sample covariance matrix and the population matrix we introduce the following quantities. Let $\{u^\nu\}_{\nu=1}^n$ and $\{\xi^\nu\}_{\nu=1}^n$ denote the realizations of the r.v. $u$ and of the noise vector $\xi$ in the given dataset $\{\mathbf{x}^\nu\}_{\nu=1}^n$. Let $\mathbf{v}_{\mathrm{PCA}}$ denote the eigenvector of the sample covariance matrix $S_n$ with largest eigenvalue $\lambda_{\mathrm{PCA}}$. Our goal is to find the relation between the finite sample values $(\lambda_{\mathrm{PCA}}, \mathbf{v}_{\mathrm{PCA}})$ and the limiting values $(\|\mathbf{v}\|^2 + \sigma^2, \mathbf{e}_1)$. Since with an exponentially small but nonzero probability $\lambda_{\mathrm{PCA}}$ may be arbitrarily small and $\langle \mathbf{v}_{\mathrm{PCA}}, \mathbf{e}_1 \rangle$ may be arbitrarily close to zero, we construct highly probable bounds on these quantities, for example, bounds that hold with probability at least $1 - \varepsilon$, where hopefully $\varepsilon \ll 1$. This is common practice both in machine learning and in concentration of measure results.

As we shall see below, the following quantities come into play in the analysis:

$$(2.6) \qquad s_u^2 = \frac{1}{n} \sum_{\nu=1}^n (u^\nu)^2, \qquad \kappa = \|\mathbf{v}\| s_u,$$

$$(2.7) \qquad \rho_j = \frac{1}{n s_u} \sum_{\nu=1}^n u^\nu \xi_j^\nu, \qquad \beta_{ij} = \frac{1}{n} \sum_{\nu=1}^n \xi_i^\nu \xi_j^\nu.$$

Loosely speaking, $s_u$ is the second moment of the variable $u$, which is close to unity for large $n$, $\kappa$ is the "signal strength" for the given dataset, the random variables $\rho_j$ capture the signal–noise interactions, and $\beta_{ij}$ are pure noise terms.

Instead of considering asymptotic results as $n \to \infty$ or as $p, n \to \infty$, in our analysis we keep $p, n$ *fixed* but view the noise level $\sigma$ as a small parameter. Further, for some of our analysis, we even keep the realizations $\{u^\nu\}$ of the random variable $u$ fixed. The following theorem provides probabilistic approximation results in terms of these quantities.

THEOREM 2.1. *Let $s_1, s_2, s_3 > 0$ and define*

$$(2.8) \qquad \varepsilon = \exp\left(-\frac{p}{2(\sqrt{5}+2)^2}\right), \qquad \varepsilon_1 = \Pr\{|N(0,1)| > s_1\},$$

$$(2.9) \qquad \varepsilon_2 = \Pr\left\{\left|\frac{\chi_{p-1}^2}{p-1} - 1\right| > \frac{s_2}{\sqrt{p-1}}\right\}, \qquad \varepsilon_3 = \Pr\{\chi_1^2 > s_3\}.$$

*Assume that $n \leq p$ and that*

$$(2.10) \qquad \kappa^2 - \frac{2\sigma s_1 \kappa}{\sqrt{n}} > \sigma^2[(1 + \sqrt{(p-1)/n})^2 + (p-1)/n];$$

*then with probability at least $1 - \varepsilon_1 - \varepsilon_2 - \varepsilon_3 - \varepsilon$,*

$$(2.11) \qquad \begin{aligned} \lambda_{\mathrm{PCA}} \geq \kappa^2 \Bigg[ 1 &- \frac{2\sigma s_1}{\kappa\sqrt{n}} + \frac{\sigma^2}{\kappa^2} \frac{p-1}{n} \frac{1 - s_2/\sqrt{p-1}}{1 + 2\sigma s_1/\kappa\sqrt{n}} \\ &- \frac{\sigma^4}{\kappa^4} \left(\frac{p-1}{n}\right)^2 \frac{(1 + s_2/\sqrt{p-1})^2}{(1 - 2\sigma s_1/\kappa\sqrt{n})^3} \Bigg] \end{aligned}$$

*and*

$$(2.12) \qquad \begin{aligned} \lambda_{\mathrm{PCA}} \leq \kappa^2 &\left(1 + \frac{2\sigma s_1}{\kappa\sqrt{n}} + \frac{\sigma^2 s_3}{\kappa^2}\right) \\ &+ \sigma\kappa\sqrt{\frac{p-1}{n}} \sqrt{1 + \frac{2\sigma s_1}{\kappa\sqrt{n}} + \frac{\sigma^2 s_3}{\kappa^2}} (1 + s_2/\sqrt{p-1}). \end{aligned}$$

*As for the leading eigenvector, with at least the same probability*

$$
\begin{aligned}
(2.13) \quad \sin\theta_{\text{PCA}} &\leq \frac{\sigma}{\kappa}\sqrt{\frac{p-1}{n}}\left(1+2\frac{\sigma s_1}{\kappa\sqrt{n}}\right)\left(1+\frac{s_2}{\sqrt{p-1}}\right) \\
&\quad + 4\sqrt{2}\frac{\sigma^2}{\kappa^2}\frac{p}{n}\frac{1}{1-(2\sigma s_1/(\kappa\sqrt{n}))-\sigma^2/\kappa^2},
\end{aligned}
$$

*where* $\sin\theta_{\text{PCA}} = \sqrt{1-\langle \mathbf{v}_{\text{PCA}}, \mathbf{e}_1\rangle^2}$.

REMARKS. Equation (2.10) can be interpreted as a condition that the signal strength is larger than the noise, since the right-hand side is a probabilistic upper bound on the norm of a noisy Wishart matrix which holds with probability at least $1-\varepsilon$; see (3.3) below. The bounds (2.11), (2.12) and (2.13) seem complicated as they involve large deviation bounds on various random variables appearing in approximations for $\lambda_{\text{PCA}}$ and $\mathbf{v}_{\text{PCA}}$. When both $p,n$ are large, $s_u$ is close to unity, so $\kappa \approx \|\mathbf{v}\|$. Assuming condition (2.10) holds, and neglecting terms $O(1/\sqrt{n})$ and $O(1/\sqrt{p})$, gives that

$$
\|\mathbf{v}\|^2 + \sigma^2\frac{p-1}{n} - \frac{\sigma^4}{\|\mathbf{v}\|^2}\left(\frac{p-1}{n}\right)^2 \lesssim \lambda_{\text{PCA}} \lesssim \|\mathbf{v}\|^2 + \sigma\|\mathbf{v}\|\sqrt{\frac{p-1}{n}}
$$

and $\sin\theta_{\text{PCA}} \lesssim (\sigma/\|\mathbf{v}\|)\sqrt{p/n} + O(\sigma^2)$. Corollary 1 below shows that the upper bound on $\sin\theta_{\text{PCA}}$ is "sharp," in the sense that $(\sigma/\|\mathbf{v}\|)\sqrt{p/n}$ is the value expected on average. Similarly, the first two terms in the lower bound for $\lambda_{\text{PCA}}$ are sharp as well.

THEOREM 2.2. *For fixed* $p,n$ *and fixed realizations* $\{u^\nu\}_{\nu=1}^n$ *and* $\{\xi^\nu\}_{\nu=1}^n$, *the largest eigenvalue and eigenvector of sample PCA are analytic functions of* $\sigma$ *inside a small interval near the origin. Inside this interval, the Taylor expansion of* $\lambda_{\text{PCA}}$ *is given by*

$$
(2.14) \quad \lambda_{\text{PCA}} = \kappa^2\left(1+2\frac{\sigma}{\kappa}\rho_1 + \frac{\sigma^2}{\kappa^2}\left(\sum_j \rho_j^2 + \beta_{11}\right) + O(\sigma^3)\right)
$$

*whereas the Taylor expansion of the corresponding eigenvector* $\mathbf{v}_{\text{PCA}}$ *is, up to a normalization constant,*

$$
(2.15) \quad \mathbf{v}_{\text{PCA}} = \mathbf{e}_1 + \frac{\sigma}{\kappa}(0,\rho_2,\ldots,\rho_p) + O(\sigma^2).
$$

Equation (2.14) shows that when $p \gg n$, even for relatively small $\sigma$, the $O(\sigma^2)$ term may be larger than the $O(\sigma)$ term, since $\rho_1/\kappa = O_P(1/(\sqrt{n}\|\mathbf{v}\|))$ whereas $1/\kappa^2\sum\rho_j^2 = O_P(p/(n\|\mathbf{v}\|^2))$.

COROLLARY 1.  *For fixed $p, n$, the mean and variance of the top eigenvalue as $\sigma \to 0$ are given by*

$$(2.16) \qquad \mathbb{E}\{\lambda_{\mathrm{PCA}}\} = \|\mathbf{v}\|^2 + \sigma^2\left(1 + \frac{p-1}{n}\right) + O(\sigma^4) + t.s.t.$$

*and*

$$(2.17) \qquad \mathrm{Var}\{\lambda_{\mathrm{PCA}}\} = \|\mathbf{v}\|^4 \frac{\mathbb{E}\{u^4\}}{n} + \frac{4\sigma^2\|\mathbf{v}\|^2\mathbb{E}\{u^2\}}{n} + O(\sigma^4) + t.s.t.$$

Here t.s.t. denotes transcendentally small terms in $\sigma$. These terms arise from the small probability of a crossover between the eigenvalue due to the signal and the largest eigenvalue due to the noise, and are of the form $A(n, p, \sigma^2)e^{-C(n,p,\|\mathbf{v}\|)/\sigma^2}$, where $A(n, p, \sigma^2)$ is analytic in $\sigma$.

Note that if $u \sim N(0, 1)$, then we obtain $\mathrm{Var}\{\lambda_{\mathrm{PCA}}\} = 2/n[\|\mathbf{v}\|^4 + 2\sigma^2\|\mathbf{v}\|^2] + O(\sigma^4)$, which recovers the asymptotic in $n$ result of Girshick for the variance of the largest eigenvalue [12].

COROLLARY 2.  *For fixed $p, n$ and fixed realizations $\{u^\nu\}_{\nu=1}^n$, the mean of $\sin\theta_{\mathrm{PCA}}$, as $\sigma \to 0$, is given by*

$$(2.18) \qquad \mathbb{E}\{\sin\theta_{\mathrm{PCA}}\} = \frac{\sigma}{\kappa\sqrt{n}} \frac{\sqrt{2}\Gamma((p-1)/2 + 1/2)}{\Gamma((p-1)/2)} + O(\sigma^2) + t.s.t.$$

The $\Gamma$ functions arise from the average of the square root of a chi-squared variable. If $p \gg 1$, then

$$\mathbb{E}\{\sin\theta_{\mathrm{PCA}}\} \approx \frac{\sigma}{\kappa\sqrt{n}}\sqrt{p}\left(1 - \frac{1}{4(p-1)} + O\left(\frac{1}{p^2}\right)\right) + O(\sigma^2) + \text{t.s.t.}$$

and

$$\mathrm{Var}\{\sin\theta_{\mathrm{PCA}}\} \approx \frac{\sigma^2}{2\kappa^2 n}(1 + O(1/p)) + O(\sigma^3) + \text{t.s.t.}$$

2.4. *Results for the joint limit $p, n \to \infty$.*  A second approach to studying PCA is to consider the joint limit $p, n \to \infty$ with $p/n = c$; see [3, 4, 9, 14, 29] and references therein. A matrix perturbation approach can also be used to prove results regarding the largest eigenvalue and corresponding eigenvector in the joint limit $p, n \to \infty$. Specifically, we obtain the following theorem.

THEOREM 2.3.  *Consider the spiked covariance model (2.4) with a single component, and assume that the random variable $u$ has finite fourth moment. Then, in the joint limit $p, n \to \infty$,*

$$(2.19) \qquad \lambda_{\mathrm{PCA}} = \begin{cases} \sigma^2\left(1 + \sqrt{\dfrac{p}{n}}\right)^2, & \text{if } n/p < \sigma^4/\|\mathbf{v}\|^4, \\[2ex] (\|\mathbf{v}\|^2 + \sigma^2)\left[1 + \dfrac{p}{n}\dfrac{\sigma^2}{\|\mathbf{v}\|^2}\right], & \text{if } n/p \geq \sigma^4/\|\mathbf{v}\|^4. \end{cases}$$

*Similarly, the dot product between the population eigenvector and the eigenvector computed by PCA also undergoes a phase transition,*

$$
(2.20) \quad R^2(p/n) = |\langle \mathbf{v}_{\text{PCA}}, \mathbf{v} \rangle|^2
$$

$$
= \begin{cases} 0, & \text{if } n/p < \sigma^4/\|\mathbf{v}\|^4, \\ \dfrac{(n\|\mathbf{v}\|^4)/(p\sigma^4) - 1}{(n\|\mathbf{v}\|^4)/(p\sigma^4) + (\|\mathbf{v}\|^2)/\sigma^2}, & \text{if } n/p \geq \sigma^4/\|\mathbf{v}\|^4. \end{cases}
$$

Equation (2.20) shows that to "learn" the direction of true largest variance, even approximately, the ratio $n/p$ must be larger than a critical threshold. This is named in the literature as *retarded learning*, or as the *phase transition phenomenon*. In the statistical physics community these results were derived using the replica method [5, 37]. In the statistics literature (2.19) was proven by Baik and Silverstein, using the Stieltjes transform [4], for the more general case of a spiked covariance model with an arbitrary finite number of independent and not necessarily Gaussian components. They showed that (with $\sigma = 1$) all eigenvalues $\alpha_j > 1 + \sqrt{p/n}$ are shifted to $\alpha_j + c\alpha_j/(\alpha_j - 1)$, and stated that "it would be interesting to have a simple heuristic argument" which shows how to obtain these pulled up values. In this section we present a matrix perturbation view of this problem, including a simple derivation of the pulled up values and some discussion as to the phenomenon for finite $p, n$. A similar approach was recently independently derived by Paul [31].

Equation (2.19) shows that for a spiked covariance model with a single component and with $\sigma = 1$, in the joint limit $p, n \to \infty$, a large enough signal eigenvalue $\alpha$ is shifted to

$$
(2.21) \quad \alpha + \frac{p}{n} \frac{\alpha}{\alpha - 1}.
$$

We now present an interesting connection between this formula and a classic result by Lawley from 1956. In [23], Lawley considered the eigenvalues of PCA for multivariate Gaussian observations whose limiting covariance matrix is diagonal with eigenvalues $\alpha_1, \ldots, \alpha_p$. Denote by $\ell_1, \ldots, \ell_p$ the noisy eigenvalues of PCA corresponding to a sample covariance matrix with a finite number of samples $n$. Then, as $n \to \infty$, with $p$ fixed

$$
(2.22) \quad \mathbb{E}\{\ell_k\} = \alpha_k + \frac{\alpha_k}{n} \sum_{i=1, i \neq k}^{p} \frac{\alpha_i}{\alpha_k - \alpha_i} + O\left(\frac{1}{n^2}\right).
$$

Applying Lawley's result to the spiked covariance model with a single significant component ($\alpha_1 = \alpha, \alpha_2 = \cdots = \alpha_p = 1$) gives

$$
\lambda_{\text{PCA}} = \alpha + \frac{p-1}{n} \frac{\alpha}{\alpha - 1} + O\left(\frac{1}{n^2}\right)
$$

whose first two terms recover the asymptotic result (2.21) of the joint limit $p$, $n \to \infty$. The remarkable point in using (2.22) is that it does not use explicit knowledge of the Marčhenko–Pastur distribution, regarding the limiting density of eigenvalues of infinitely large random matrices. We remark that although the first two terms in Lawley's expansion provide the asymptotic result as $p, n \to \infty$, this is not due to the higher-order terms all vanishing individually. Rather, in the joint limit $p, n \to \infty$ they all miraculously cancel each other. Yet, based on this observation, we propose the following two theorems.

THEOREM 2.4. *Consider a more general Gaussian spiked covariance model with $k$ large components with fixed variances $\alpha_1, \ldots, \alpha_k$ and with the $p - k$ remaining components each having a random variance sampled independently from a density $h(\alpha)$ with compact support in the interval $[0, \alpha_c]$ ($h(\alpha) = 0$ for $\alpha > \alpha_c$). Then, in the joint limit $p, n \to \infty$ and for large enough values $\alpha_j$, the first $k$ largest limiting eigenvalues of PCA converge to*

$$(2.23) \qquad \lambda_j = T(\alpha_j) = \alpha_j + \frac{p}{n} \alpha_j \int_0^{\alpha_c} \frac{\rho}{\alpha_j - \rho} h(\rho) \, d\rho.$$

THEOREM 2.5. *Let $\{\alpha_j\}$ denote an infinite sequence of i.i.d. random variables from a density $h(\alpha)$ with compact support. Let $\mathbf{x} = (x_1, \ldots, x_p)$ be a $p$-dimensional vector composed of $p$ independent Gaussian random variables, where each $x_j$ has variance $\alpha_j$. Let $C(n, p)$ denote the $p \times p$ empirical covariance matrix computed from $n$ independent samples $\{\mathbf{x}_i\}_{i=1}^n$ from this model. Then, in the joint limit $p, n \to \infty$, $p/n = c$, the spectral norm of $C$ is equal to*

$$(2.24) \qquad \lim_{p, n \to \infty} \|C\| = \alpha^* + c\alpha^* \int \frac{\rho}{\alpha^* - \rho} h(\rho) \, d\rho$$

*where $\alpha^*$ is the maximal point at which*

$$(2.25) \qquad \frac{dT(\alpha)}{d\alpha}\Big|_{\alpha^*} = 0.$$

A motivation for the model considered in these theorems is a setting where high-dimensional observations are of the type "signal plus noise," but where the noise is heteroscedastic and possibly correlated. Thus, in a suitable basis, different noise components have different variances, and we only know some statistical properties about the noise, such as the distribution of these variances. Theorem 2.4 is a generalization of (2.21) that holds for the standard spiked covariance model. Theorem 2.5 follows immediately from Theorem 2.4 according to the following reasoning: In this modified model there is also a similar phase transition phenomenon. If the original variance of the signal, $\alpha$, is larger than a critical value $\alpha^*(h, p, n)$, then it will be pulled up from the noise in the limit $p, n \to \infty$. Further, for all $\alpha > \alpha^*$, this pulled up value is monotonic in $\alpha$. From Theorem 2.4, the critical value $\alpha^*$ satisfies (2.25), and at that point according to our matrix perturbation analysis, the value $T(\alpha^*)$ is equal to the spectral norm $C$—the covariance matrix of the noise. We remark that a formula similar to (2.24) was recently derived by El Karoui [11], who also studied the finite $p, n$ fluctuations around this mean.

COROLLARY 3. *Consider the general spiked covariance model as in Theorem 2.4, and assume $c = p/n \gg 1$. Let*

$$\mu_1 = \int \rho h(\rho)\, d\rho, \qquad \mu_2^2 = \int (\rho - \mu_1)^2 h(\rho)\, d\rho.$$

*Then, in the joint limit $p, n \to \infty$, the norm of a pure noise matrix is approximately*

$$(2.26) \qquad \lambda(\alpha^*) = \mu_1\left(c + 2\sqrt{c}\sqrt{1 + \frac{\mu_2^2}{\mu_1^2}} + O(1)\right).$$

*Further, the phase transition phenomenon for the pulled up eigenvalues occurs at*

$$(2.27) \qquad \alpha^* = \mu_1\left(\sqrt{c}\sqrt{1 + \frac{\mu_2^2}{\mu_1^2}} + 1 + 4\frac{\mu_2^2/\mu_1^2}{1 + 2\mu_2^2/\mu_1^2} + O(1/\sqrt{c})\right).$$

Equation (2.26) may be useful for inference on the number of components in a general spiked covariance model, given that the first two moments $\mu_1, \mu_2$ of the density $h(\rho)$ of the noise are either known a priori or estimated from the data.

## 3. Proof of Theorem 2.1.

To prove Theorem 2.1 we shall use the following three lemmas:

LEMMA 1. *Let $A, B$ be $p \times p$ Hermitian matrices. Let $\{\lambda_i\}_{i=1}^p$ denote the eigenvalues of $A$ sorted in decreasing order with corresponding eigenvectors $\mathbf{v}_i$. Let $P_i$ denote the projection into the orthogonal subspace of $\mathbf{v}_i$, $P_i\mathbf{v} = \mathbf{v} - \langle \mathbf{v}, \mathbf{v}_i \rangle \mathbf{v}_i$. If $\lambda_1$ has multiplicity 1 and $\|B\| < \lambda_1 + \langle \mathbf{v}_1, B\mathbf{v}_1 \rangle - \lambda_2$, then the largest eigenvalue of $A + B$ satisfies the bounds*

$$(3.1) \qquad \langle \mathbf{v}_1, B\mathbf{v}_1 \rangle \leq \lambda_1(A + B) - \lambda_1 \leq \langle \mathbf{v}_1, B\mathbf{v}_1 \rangle + \|P_1 B\mathbf{v}_1\|.$$

LEMMA 2. *Let $X$ denote an $n \times p$ matrix with entries $X_{ij}$ all i.i.d. $N(0, 1)$ Gaussian variables, and let $W = X^T X/n$ be the corresponding scaled Wishart matrix. Define*

$$(3.2) \qquad \varepsilon = \exp\left(-\frac{p}{2(\sqrt{5} + 2)^2}\right),$$

*then for $n \leq p$ with probability at least $1 - \varepsilon$,*

$$\|W\| \leq (1 + \sqrt{p/n})^2 + p/n \tag{3.3}$$

*and*

$$\|W - I_p\| \leq 4\frac{p}{n}. \tag{3.4}$$

LEMMA 3. *Let $A$ be a $p \times p$ Hermitian matrix and let $B$ be a Hermitian perturbation. Let $(\lambda, \mathbf{v})$ be the eigenvalue/vector pair of $A + B$ corresponding to $(\lambda_i, \mathbf{v}_i)$ of $A$, and let $\delta = \min_{j \neq i} |\lambda - \lambda_j|$, where $\{\lambda_j\}_{j=1}^p$ are all the eigenvalues of $A$; then*

$$\sin \theta(\mathbf{v}, \mathbf{v}_i) \leq \frac{\|B\|}{\delta}. \tag{3.5}$$

Lemma 1 follows from classical results in matrix or operator perturbation theory. According to [30], Theorem 4.5.1 (see also [13], Theorem 6.3.14), for each scalar $\mu$, vector $\mathbf{x}$ and Hermitian matrix $A$, there exists an eigenvalue $\lambda$ of $A$ such that

$$|\lambda - \mu| \leq \|A\mathbf{x} - \mu\mathbf{x}\|/\|\mathbf{x}\|.$$

Applying this theorem to the matrix $A + B$, with $\mathbf{x} = \mathbf{v}_i$ a normalized eigenvector of $A$ and with $\mu = \lambda_i + \langle \mathbf{v}_i, B\mathbf{v}_i \rangle$, gives that

$$|\lambda(A + B) - \lambda_i(A) - \langle \mathbf{v}_i, B\mathbf{v}_i \rangle| \leq \|P_i B\mathbf{v}_i\|.$$

The condition of Lemma 1 ensures that the largest eigenvalue of $A$ is the one closest to the largest eigenvalue of $A + B$, for example, $i = 1$. For the analysis of a spiked covariance model with more than one component, similar statements can be made for interior eigenvalues as well [16]. Lemma 2 follows from Theorem II.13 of Szarek and Davidson [7] and is proved in the Appendix. Lemma 3 is known as the $\sin \theta$ theorem of Davis and Kahan [8]; see also [30], Theorem 11.7.1.

PROOF OF THEOREM 2.1.    Let $\{\mathbf{x}^\nu\}_{\nu=1}^n$ be $n$ i.i.d. observations from the one-component model (2.4). We decompose the corresponding sample covariance matrix as follows:

$$S_n = \begin{pmatrix} \kappa^2 & 0 & \cdots & 0 \\ 0 & 0 & & 0 \\ \vdots & & & \vdots \\ 0 & 0 & \cdots & 0 \end{pmatrix} + \sigma\kappa \begin{pmatrix} 2\rho_1 & \rho_2 & \cdots & \rho_p \\ \rho_2 & 0 & & 0 \\ \vdots & & 0 & \vdots \\ \rho_p & 0 & & 0 \end{pmatrix}$$

$$(3.6) \qquad + \sigma^2 \begin{pmatrix} \beta_{1,1} & \beta_{1,2} & \cdots & \beta_{1,p} \\ \beta_{2,1} & \beta_{2,2} & & \vdots \\ \vdots & & \ddots & \\ \beta_{p,1} & \cdots & & \beta_{p,p} \end{pmatrix}$$

$$= \mathcal{L}_0 + \sigma \mathcal{L}_1 + \sigma^2 \mathcal{L}_2,$$

where $\kappa, \rho_j$ and $\beta_{ij}$ are defined above in (2.6) and (2.7). Note that conditional on the realizations $\{u^\nu\}_{\nu=1}^n$ of the random variable $u$ kept fixed, $\rho_i = \frac{1}{\sqrt{n}} \eta_i$, where $\eta_i$ are all i.i.d. $N(0,1)$, $\beta_{ii} = \chi_n^2/n$ and that $\mathcal{L}_2$ is a Wishart noise matrix. The matrix $\mathcal{L}_0$ can be thought of as the "signal," the matrix $\mathcal{L}_1$ as the signal–noise interactions, whereas $\mathcal{L}_2$ contains pure noise.

To derive the lower bound (2.11), we view the matrix $\sigma^2 \mathcal{L}_2$ as a *perturbation* of the matrix $\mathcal{L}_0 + \sigma \mathcal{L}_1$. Since the matrix $\mathcal{L}_2$ is nonnegative (it is a scaled Wishart matrix), it follows that

$$\lambda_{\mathrm{PCA}} \geq \|\mathcal{L}_0 + \sigma \mathcal{L}_1\|.$$

The matrix $\mathcal{L}_0 + \sigma \mathcal{L}_1$ has rank 2 with the following two nonzero eigenvalues:

$$(3.7) \qquad \lambda_{\pm} = \frac{(\kappa^2 + 2\sigma\kappa\rho_1) \pm \sqrt{(\kappa^2 + 2\sigma\kappa\rho_1)^2 + 4(\sigma\kappa)^2 \sum_{j \geq 2} \rho_j^2}}{2},$$

where $\lambda_+$ is positive and $\lambda_-$ is negative. Since conditional on the realizations $u^\nu$ fixed, the random variables $\rho_j$ are independent Gaussians, we define

$$(3.8) \qquad \sum_{j \geq 2} \rho_j^2 = \frac{T_1}{n},$$

where $T_1 \sim \chi_{p-1}^2$.

Using the inequality $\sqrt{1+x} \geq (1 + x/2 - x^2/8)$ in (3.7) gives

$$\lambda_+ \geq (\kappa^2 + 2\sigma\kappa\rho_1) \left[ 1 + \frac{\sigma^2 \kappa^2}{(\kappa^2 + 2\sigma\kappa\rho_1)^2} \frac{T_1}{n} - \frac{(\sigma\kappa)^4}{(\kappa^2 + 2\sigma\kappa\rho_1)^4} \frac{T_1^2}{n^2} \right].$$

Therefore, with probability at least $1 - \varepsilon_1 - \varepsilon_2$

$$\lambda_{\mathrm{PCA}} \geq \lambda_+ \geq \kappa^2 \left( 1 - \frac{2\sigma s_1}{\kappa\sqrt{n}} \right) \left[ 1 + \frac{\sigma^2}{\kappa^2} \frac{p-1}{n} \frac{1 - s_2/\sqrt{p-1}}{1 + (2\sigma s_1)/(\kappa\sqrt{n})} \right.$$

$$\left. - \frac{\sigma^4}{\kappa^4} \frac{(p-1)^2}{n^2} \frac{(1 + s_2/\sqrt{p-1})^2}{(1 - (2\sigma s_1)/(\kappa\sqrt{n}))^3} \right],$$

which proves the lower bound (2.11).

To prove the upper bound on $\lambda_{\mathrm{PCA}}$, we use Lemma 1 with the matrix $\sigma \mathcal{L}_1 + \sigma^2 \mathcal{L}_2$ as a perturbation of $\mathcal{L}_0$. Choosing $\mathbf{e}_1$ as an eigenvector of $\mathcal{L}_0$ and applying the lemma gives that

$$\lambda_{\text{PCA}} \leq \kappa^2 + 2\sigma\kappa\rho_1 + \sigma^2\beta_{11} + \sigma\sqrt{\sum_{j \geq 2}(\kappa\rho_j + \sigma\beta_{1j})^2}.$$

Conditional on fixed realizations $u^\nu$ and on fixed noise realizations $\xi_1^\nu$ in the first component of the data, we have that

$$\kappa\rho_j + \sigma\beta_{1j} = \frac{1}{\sqrt{n}}\sqrt{\kappa^2 + 2\sigma\kappa\rho_1 + \sigma^2\beta_{11}}\,\eta_j,$$

where the $\eta_j$ are independent standard Gaussian variables. Therefore,

$$\lambda_{\text{PCA}} \leq (\kappa^2 + 2\sigma\kappa\rho_1 + \sigma^2\beta_{11}) + \sigma\sqrt{\kappa^2 + 2\sigma\kappa\rho_1 + \sigma^2\beta_{11}}\sqrt{\frac{T}{n}},$$

where the random variable $T \sim \chi_{p-1}^2$, and is independent of $\rho_1$ and $\beta_{11}$. Therefore, with probability at least $1 - \varepsilon - \varepsilon_1 - \varepsilon_2 - \varepsilon_3$ the bound of (2.12) follows.

To prove a bound on the eigenvector $\mathbf{v}_{\text{PCA}}$, we start from the eigenvector $\mathbf{v}_+$ corresponding to $\lambda_+$, and given by

$$(3.9) \qquad \mathbf{v}_+ = \frac{1}{Z}\left[\mathbf{e}_1 + \frac{\sigma\kappa}{\lambda_+}(0, \rho_2, \ldots, \rho_p)^T\right],$$

where $Z = \sqrt{1 + \sigma^2\kappa^2 T_1/n\lambda_+^2}$ is a normalization constant such that $\|\mathbf{v}_+\| = 1$.

Simple algebraic manipulations and the triangle inequality give

$$\sin\theta_{\text{PCA}} = \sqrt{1 - \langle\mathbf{v}_{\text{PCA}}, \mathbf{e}_1\rangle^2} = \sqrt{1 + |\langle\mathbf{v}_{\text{PCA}}, \mathbf{e}_1\rangle|}\sqrt{1 - |\langle\mathbf{v}_{\text{PCA}}, \mathbf{e}_1\rangle|}$$

$$(3.10) \qquad = \sqrt{1 + |\langle\mathbf{v}_{\text{PCA}}, \mathbf{e}_1\rangle|}\frac{\|\mathbf{v}_{\text{PCA}} - \mathbf{e}_1\|}{\sqrt{2}}$$

$$\leq \|\mathbf{v}_{\text{PCA}} - \mathbf{e}_1\| \leq \|\mathbf{v}_{\text{PCA}} - \mathbf{v}_+\| + \|\mathbf{v}_+ - \mathbf{e}_1\|.$$

From (3.9), a bound on the second term in (3.10) is

$$(3.11) \qquad \|\mathbf{v}_+ - \mathbf{e}_1\| = \sqrt{2}\sqrt{1 - \frac{1}{Z}} \leq \frac{\sigma\kappa}{\lambda_+}\sqrt{\frac{T_1}{n}}.$$

For the first term in (3.10), applying the $\sin\theta$ theorem (Lemma 3 above) with the matrix $\sigma^2(W - I_p) = \sigma^2(\mathcal{L}_2 - I_p)$ as the perturbation of $\mathcal{L}_0 + \sigma\mathcal{L}_1 + \sigma^2 I_p$ gives

$$(3.12) \qquad \|\mathbf{v}_{\text{PCA}} - \mathbf{v}_+\| = \sqrt{2}\sqrt{1 - |\langle\mathbf{v}_{\text{PCA}}, \mathbf{v}_+\rangle|} \leq \sqrt{2}\sin\theta(\mathbf{v}_{\text{PCA}}, \mathbf{v}_+)$$

$$\leq \sqrt{2}\sigma^2\frac{\|W - I_p\|}{\delta},$$

where $\delta = \min_{j \neq 1} |\lambda_{\mathrm{PCA}} - \lambda_j(\mathcal{L}_0 + \sigma\mathcal{L}_1 + \sigma^2 I_p)|$. Therefore, combining (3.11) and (3.12),

$$(3.13) \qquad \sin\theta_{\mathrm{PCA}} \leq \frac{\sigma\kappa}{\lambda_+}\sqrt{\frac{T_1}{n}} + \sigma^2\sqrt{2}\frac{\|W - I_p\|}{\delta}.$$

To conclude the proof we apply almost sure bounds for $\|W - I_p\|$ and $T_1$ from above and $\delta$ and $\lambda_+$ from below. Bounds on $T_1$ and $\lambda_+$ are identical to those described in the proof of (2.11) of the theorem. To bound $\|W - I_p\|$ from above, we use (3.4) of Lemma 2, which states that with probability at least $1 - \varepsilon$

$$(3.14) \qquad \|W - I_p\| \leq 4\frac{p}{n}.$$

As for a bound on $\delta$, from the first part of the proof, if $\kappa^2 + 2\sigma\kappa\rho_1 > \sigma^2[(1 + \sqrt{p/n})^2 + p/n]$,

$$\lambda_{\mathrm{PCA}} \geq \lambda_+ \geq \kappa^2 + 2\sigma\kappa\rho_1.$$

Furthermore, the eigenvalues of the matrix $\mathcal{L}_0 + \mathcal{L}_1 + \sigma^2 I$ are $(\lambda_+) + \sigma^2, \sigma^2$ or $(\lambda_-) + \sigma^2$, with $\lambda_\pm$ given by (3.7). Therefore,

$$\delta = \min_{j \neq 1} |\lambda_{\mathrm{PCA}} - \lambda_j(\mathcal{L}_0 + \mathcal{L}_1 + \sigma^2 I)| = \lambda_{\mathrm{PCA}} - \sigma^2 \geq \kappa^2 + 2\sigma\kappa\rho_1 - \sigma^2$$

and with probability at least $1 - \varepsilon_1$,

$$(3.15) \qquad \delta \geq \kappa^2 - 2s_1\frac{\sigma\kappa}{\sqrt{n}} - \sigma^2.$$

Combining (3.13), (3.14) and (3.15) concludes the proof. $\square$

**4. Proof of Theorem 2.2.** We now explore the leading order terms in $\sigma$ of the explicit dependence of $\lambda_{\mathrm{PCA}}$ and $\mathbf{v}_{\mathrm{PCA}}$ on $p, n$, and on the specific signal and noise realizations $\{u^\nu\}_{\nu=1}^n$ and $\{\xi^\nu\}_{\nu=1}^n$. To this end, we view $\sigma$ as a small parameter and consider the Taylor expansion of $\lambda_{\mathrm{PCA}}$ and $\mathbf{v}_{\mathrm{PCA}}$ as $\sigma \to 0$. By definition, the largest eigenvalue $\lambda_{\mathrm{PCA}}$ is the largest root of the characteristic polynomial of sample covariance matrix $S_n$. For $\sigma = 0$ this eigenvalue is a simple root with multiplicity 1. Therefore, given a finite dataset $\{\mathbf{x}^\nu\}_{\nu=1}^n$ with $s_u > 0$, the largest eigenvalue $\lambda(\sigma)$, when viewed as a function of noise level, is an *analytic* function of $\sigma$ in the complex plane for small enough $\sigma$. This statement follows from the representation of the empirical covariance matrix as $S_n = \mathcal{L}_0 + \sigma\mathcal{L}_1 + \sigma^2\mathcal{L}_2$, (3.6), with all matrices being symmetric, together with standard results regarding perturbation theory for linear operators; see, for example, Kato [21], Chapter 2, Theorem 6.1. Moreover, the radius of convergence of a Taylor series of $\lambda(\sigma)$ is the largest complex $\sigma$ for which $\lambda_{\mathrm{PCA}}(\sigma) > \lambda_2(\sigma)$ where $\lambda_2(\sigma)$ is the second largest eigenvalue of the noisy covariance matrix. Note that the location of the crossover depends on the specific signal and noise realizations. Finally, since the matrix $S_n$ is symmetric $\lambda(\sigma)$ is real when $\sigma$ is real.

Therefore, for small enough $\sigma$ we can expand both the top eigenvalue and its corresponding eigenvector as a regular power series in $\sigma$:

$$\mathbf{v}_{\mathrm{PCA}} = \mathbf{v}_0 + \sigma \mathbf{v}_1 + \sigma^2 \mathbf{v}_2 + \cdots,$$

$$\lambda_{\mathrm{PCA}} = \lambda_0 + \sigma \lambda_1 + \sigma^2 \lambda_2 + \cdots.$$

We insert these expansions into (3.6) and equate terms with equal powers of $\sigma$. The first few equations read

$$\mathcal{L}_0 \mathbf{v}_0 = \lambda_0 \mathbf{v}_0,$$

$$\mathcal{L}_0 \mathbf{v}_1 + \mathcal{L}_1 \mathbf{v}_0 = \lambda_0 \mathbf{v}_1 + \lambda_1 \mathbf{v}_0,$$

$$\mathcal{L}_0 \mathbf{v}_2 + \mathcal{L}_1 \mathbf{v}_1 + \mathcal{L}_2 \mathbf{v}_0 = \lambda_0 \mathbf{v}_2 + \lambda_1 \mathbf{v}_1 + \lambda_2 \mathbf{v}_0.$$

Iteratively solving these equations gives

$$(4.1) \qquad \lambda = \kappa^2 + 2\sigma\kappa\rho_1 + \sigma^2 \left( \sum_{j \geq 2} \rho_j^2 + \beta_{11} \right) + O(\sigma^3)$$

and

$$(4.2) \qquad \begin{aligned} \mathbf{v} &= \mathbf{e}_1 + \frac{\sigma}{\kappa}(0, \rho_2, \ldots, \rho_p) \\ &\quad + \frac{\sigma^2}{\kappa^2}[(0, \beta_{12}, \ldots, \beta_{1p}) - 2\rho_1(0, \rho_2, \ldots, \rho_p)] + O\left(\frac{\sigma^3}{\kappa^3}\right). \end{aligned}$$

Note that up to order $O(\sigma^2)$, the eigenvalue $\lambda_{\mathrm{PCA}}$ and the corresponding eigenvector depend only on the first row of the noisy matrix, for example, only on the interaction between signal and noise.

Equations (2.16), (2.17) and (2.18) follow by taking expectations on expressions (4.1) and (4.2) for $\lambda_{\mathrm{PCA}}$ and $\mathbf{v}_{\mathrm{PCA}}$, respectively, and retaining only the leading terms in $\sigma$. However, an important remark is that (4.1) and (4.2) are the Taylor expansions of the eigenvalue and eigenvector that are analytic in $\sigma$ and correspond to $\kappa^2$ and $\mathbf{e}_1$ when $\sigma = 0$. As such, these are not necessarily expansions of $\lambda_{\mathrm{PCA}}$ and $\mathbf{v}_{\mathrm{PCA}}$—the actual largest eigenvalue and corresponding eigenvector of the sample covariance matrix. This is because for finite $\sigma > 0$ there is a nonzero probability that the largest eigenvalue is one due to noise and not the one described by (4.1). From Lemma 2, the probability of such a crossover, between the eigenvalue due to noise and the eigenvalue due to the signal, can be bounded by an expression of the form $A(n, p) \exp(-C(n, p)/\sigma^2)$. Therefore, by taking expectations of (4.1), we introduce *transcendentally small error terms* in $\sigma$.

**5. Proof of Theorem 2.3: the phase transition phenomenon.**

5.1. *A simple heuristic for the location of the phase transition.* First, we present a simple heuristic explanation for the phase transition phenomenon, but for fixed $p, n$, as a function of noise level $\sigma$. Obviously, for fixed $p, n$ there is no deterministic phase transition at a fixed constant $c = p/n$, only an increasing probability for losing the relation between the direction of maximal variance and the limiting vector $\mathbf{e}_1$. From the analysis of Section 3, this occurs when the largest eigenvalue of the sample covariance matrix is of the same order of magnitude as that of the noise matrix $E$, $\lambda_{\mathrm{PCA}} \sim \|E\|$. From (3.4), this occurs roughly when

$$\kappa^2 + \sigma^2 + \sigma^2 \frac{p}{n} = \sigma^2 (1 + \sqrt{p/n})^2.$$

This gives

$$\frac{p}{n} = \frac{1}{4} \frac{\kappa^4}{\sigma^4},$$

which up to a multiplicative constant has the same functional dependence on the parameters $p, n, \sigma$ as the true location for the phase transition in (2.19) and (2.20).

5.2. *An exact analysis of the phase transition.* We now present a simple linear-algebra based derivation of the exact pulled up value for a spiked population model with a single component. A similar though more complicated analysis applies for the general $k$-component model. For simplicity, we perform our analysis for $p/n = c = 1$, and without loss of generality assume $\sigma = 1$.

To this end, we decompose the $p \times p$ sample covariance matrix computed from $n$ samples as follows:

$$S_n = \begin{pmatrix} \kappa^2 + 2\kappa\rho_{11} + \beta_{11} & b_2 & \cdots & b_p \\ b_2 & 0 & \cdots & 0 \\ \vdots & \vdots & \ddots & \vdots \\ b_p & 0 & \cdots & 0 \end{pmatrix} + \begin{pmatrix} 0 & 0 & \cdots & 0 \\ 0 & \beta_{22} & \cdots & \beta_{2p} \\ \vdots & \beta_{32} & \ddots & \beta_{3p} \\ 0 & \beta_{p2} & \cdots & \beta_{pp} \end{pmatrix},$$

where $b_j = \kappa\rho_j + \beta_{1j}$. Note that the second matrix, which is the minor of the covariance matrix obtained by deleting the first row and column, is just a $(p-1) \times (p-1)$ scaled Wishart matrix. Let $\lambda_2, \ldots, \lambda_p$ denote its eigenvalues and let $V_{p \times p}$ be the matrix of its eigenvectors padded with zeros in the first coordinate. We perform a change of basis whereby this matrix becomes diagonal. Then the full covariance matrix takes the following form:

$$(5.1) \qquad VS_nV^T = \begin{pmatrix} \kappa^2 + 2\kappa\rho_1 + \beta_{11} & \tilde{b}_2 & \cdots & \tilde{b}_p \\ \tilde{b}_2 & \lambda_2 & \cdots & 0 \\ \vdots & \vdots & \ddots & \vdots \\ \tilde{b}_p & 0 & \cdots & \lambda_p \end{pmatrix},$$

where $\tilde{b}_{1j}$ are the entries of the first row and column in the new basis. The specific form (5.1) is known as an *arrowhead matrix* [28]. Assuming that $\tilde{b}_j \neq 0$ for all $j$ (an event with probability 1), the eigenvalues of this matrix satisfy the following *secular equation*:

$$(5.2) \qquad f(\lambda) = (\lambda - \kappa^2 - 2\kappa\rho_1 - \beta_{11}) - \sum_{j=2}^{p} \frac{\tilde{b}_j^2}{\lambda - \lambda_j} = 0.$$

Recall that $\tilde{b}_j = \tilde{\rho}_j + \tilde{\beta}_{1j}$ are the entries of the first row and column in the new basis. They are given explicitly as

$$\tilde{\rho}_j = \frac{1}{n s_y} \sum_{\nu=1}^{n} u_j^\nu \tilde{\xi}_j^\nu, \qquad \tilde{\beta}_{1j} = \frac{1}{n} \sum_{\nu=1}^{n} \tilde{\xi}_1^\nu \tilde{\xi}_j^\nu,$$

where $\tilde{\xi}^\nu$ are the noise vectors in the rotated basis in which the $(p-1) \times (p-1)$ submatrix of noise covariances is diagonal with eigenvalues $\lambda_j$. Therefore, conditional on $\tilde{\xi}_j$ having variance $\lambda_j$, the quantity $\tilde{b}_j$ is $N(0, \lambda_j(\kappa^2 + 2\kappa\rho_1 + \beta_{11})/n)$. Therefore,

$$\sum_{j=2}^{p} \frac{\tilde{b}_j^2}{\lambda - \lambda_j} = \frac{p-1}{n}(\kappa^2 + 2\kappa\rho + \beta_{11})\frac{1}{p-1}\sum_{j=2}^{p} \frac{\lambda_j \eta_j^2}{\lambda - \lambda_j},$$

where $\eta_j$ are all i.i.d. $N(0,1)$ and independent of $\lambda_j$. Furthermore, in the limit $p, n \to \infty, p/n = c$, the distribution of eigenvalues $\lambda_j$ of the submatrix converges to the Marčenko–Pastur distribution [24],

$$f_{\mathrm{MP}}(x) = \frac{1}{2\pi x c}\sqrt{(b-x)(x-a)}, \qquad x \in [a,b],$$

where $a = (1 - \sqrt{c})^2, b = (1 + \sqrt{c})^2$. In addition, as $n, p \to \infty$, $\kappa^2 \to \|\mathbf{v}\|^2$, $\rho_1 = O_P(1/\sqrt{n}) \to 0$ and $\beta_{11} = \chi_n^2/n \to 1$, all with probability 1. Therefore, as $p, n \to \infty$, the sum in (5.2) converges with probability 1 to the following integral:

$$(5.3) \qquad \lim_{p,n\to\infty} \sum_{j=2}^{p} \frac{\tilde{b}_j^2}{\lambda - \lambda_j} = (\|\mathbf{v}\|^2 + 1)\frac{p-1}{n}\int_a^b f_{\mathrm{MP}}(x)\frac{x}{\lambda - x}\,dx.$$

This integral is a linear functional of the Marčenko–Pastur distribution, with some similarity to its Stieltjes transform. We remark that the Stieltjes transform was used extensively in deriving results regarding the eigenvalues of random matrices [24, 33].

This integral can be evaluated explicitly. For example, for $c = 1$,

$$(5.4) \qquad \int_0^4 \frac{1}{2\pi}\sqrt{(4-x)x}\frac{1}{\lambda - x}\,dx = \frac{1}{2}[\lambda - 2 - \sqrt{\lambda(\lambda - 4)}].$$

Thus, the largest eigenvalue satisfies a quadratic equation in $\lambda$, whose solution is

$$\lambda(\alpha) = \alpha + c\frac{\alpha}{\alpha - 1}$$

with $\alpha = \|\mathbf{v}\|^2 + 1$. However, this solution is indeed the largest eigenvalue only if $\lambda(\alpha) > (1 + \sqrt{c})^2$. This recovers the pulled up value and the exact location of the phase transition, (2.19).

To prove (2.20) for the eigenvector overlap, we note that eigenvectors of arrowhead matrices have also a closed form expression [28]. Let $\lambda$ be an eigenvalue of the arrowhead matrix (5.1); then up to normalization, the corresponding eigenvector is given by

$$(5.5) \qquad \mathbf{v} = \left(1, \frac{\tilde{b}_2}{\lambda - \lambda_2}, \ldots, \frac{\tilde{b}_p}{\lambda - \lambda_p}\right).$$

Therefore,

$$R^2 = \frac{\langle \mathbf{v}, \mathbf{e}_1 \rangle^2}{\|\mathbf{v}\|^2} = \frac{1}{1 + \sum_{j \geq 2} \tilde{b}_j^2/(\lambda - \lambda_j)^2}.$$

In the joint limit $p, n \to \infty$, similar to the analysis above, the sum in the denominator converges with probability 1 to an analogous integral

$$(5.6) \qquad \lim_{p,n\to\infty, p/n=c} R^2 = \frac{1}{1 + p/n \int \alpha\mu/(\lambda - \mu)^2 f_{\mathrm{MP}}(\mu)\,d\mu}.$$

Evaluation of the integral gives the asymptotic overlap, (2.20).

5.3. *A classical result of Lawley and two theorems.* While Theorems 2.4 and 2.5 are motivated by the classical result of Lawley, (2.22), their proof relies on results from random matrix theory regarding the limiting empirical density of eigenvalues of sample covariance matrices in the joint limit $p, n \to \infty$. Before proving these theorems, we first briefly review the results required for our proofs.

The key required quantity is the Stieltjes transform of a probability density $h(t)$ defined as

$$(5.7) \qquad m_h(z) = \int \frac{h(t)}{t - z}\,dt \qquad \forall z \in \mathbb{C}^+.$$

Let $S_n = 1/n Z^T Z$ be the $p \times p$ sample covariance matrix of $n$ observations $\mathbf{z}_i \in \mathbb{R}^p$ from the model described in Theorem 2.5, and denote by $F_n$ the empirical distribution function of the eigenvalues of $S_n$,

$$(5.8) \qquad F_n(t) = \frac{\{\#\mu_j < t\}}{p}.$$

As proven in [33], in the limit $p, n \to \infty$, $F_n$ converges with probability 1 to a limiting distribution $F$ with no eigenvalues of $S_n$ outside the support of this limiting distribution. Although the explicit form of $F$ can be computed only in a handful of simple cases, its Stieltjes transform satisfies the equation

$$(5.9) \qquad m(z) = \int \frac{h(t)}{t(1 - c - m(z)z) - z} \, dt.$$

One can also consider a different matrix, $\bar{S}_n = 1/n ZZ^T$ of size $n \times n$. Since the matrices $S_n$ and $\bar{S}_n$ have the same nonzero eigenvalues and differ by $|p - n|$ zero eigenvalues, their respective empirical distribution functions $F$ and $\bar{F}$ are related as follows:

$$\bar{F} = (1 - c) I_{[0, \infty)} + cF.$$

Due to linearity of the Stieltjes transform,

$$(5.10) \qquad \bar{m}(z) = -\frac{1 - c}{z} + cm(z).$$

Obviously, when $c = 1$, $\bar{m}(z) = m(z)$.

The last result of interest is an inverse equation for $\bar{m}(z)$, which reads

$$(5.11) \qquad z(\bar{m}) = -\frac{1}{\bar{m}} + c \int \frac{t}{1 + t\bar{m}} h(t) \, dt.$$

PROOF OF THEOREM 2.4. For simplicity, we consider a spiked covariance model with a single spike, assumed in the direction $\mathbf{e}_1$. Let $\alpha_1$ be fixed and sufficiently large, and let $\{\alpha_j\}_{j \geq 2}$ denote a sequence of i.i.d. realizations sampled from a density $h(\alpha)$. Consider $n$ i.i.d. random vectors $\{\mathbf{x}_\nu\}_{\nu=1}^n$ from a model

$$\mathbf{x} = \sqrt{\alpha_1} y_1 \mathbf{e}_1 + \sum_{j=2}^p \sqrt{\alpha_j} y_j \mathbf{e}_j,$$

where $y_j$ are all i.i.d. Gaussian $N(0, 1)$ random variables. We view the direction $\mathbf{e}_1$ as the signal direction with all other directions as noise, and write the corresponding sample covariance matrix as follows:

$$S_n = \begin{pmatrix} \kappa^2 & b_1 & \cdots & b_p \\ b_1 & & & \\ \vdots & & C_n & \\ b_p & & & \end{pmatrix},$$

where

$$\kappa^2 = \frac{1}{n} \sum_{\nu=1}^n (x_1^\nu)^2, \qquad b_j = \frac{1}{n} \sum_{\nu=1}^n x_1^\nu x_j^\nu$$

and $C_n$ is the $(p-1) \times (p-1)$ sample covariance matrix of the pure noise components.

Let $V_0$ be a $(p-1) \times (p-1)$ matrix that diagonalizes the pure noise matrix $C_n$, and consider the $p \times p$ unitary matrix

$$V = \left( \begin{array}{c|c} 1 & \\ \hline & V_0 \end{array} \right).$$

In the basis defined by $V$, the sample covariance matrix takes the form

$$(5.12) \qquad V S_n V^{-1} = \begin{pmatrix} \kappa^2 & \tilde{b}_1 & \tilde{b}_2 & \cdots & \tilde{b}_p \\ \tilde{b}_1 & \mu_1 & & & \\ \tilde{b}_2 & & \mu_2 & & \\ \vdots & & & \ddots & \\ \tilde{b}_p & & & & \mu_p \end{pmatrix},$$

where $\tilde{b}_j$ is the projection of the vector $\mathbf{b}$ on the $j$th basis vector of the matrix $V$,

$$\tilde{b}_j = \mathbf{b}^T V_j = \frac{1}{n} \sum_{\nu=1}^{n} x_1^{\nu} \tilde{x}_j^{\nu}.$$

As in (5.1), the matrix (5.12) has the form of a symmetric arrowhead matrix, and assuming all $\tilde{b}_j \neq 0$, its eigenvalues are solutions of

$$(5.13) \qquad \lambda - \kappa^2 = \sum_{j=1}^{p} \frac{\tilde{b}_j^2}{\lambda - \mu_j}.$$

We now consider the joint limit $p, n \to \infty, p/n = c$. Similarly to the analysis of Section 5.2, the sum in (5.13) converges with probability 1 to

$$(5.14) \qquad c\alpha_1 \int \frac{\mu}{\lambda - \mu} \, dF(\mu),$$

where $F(\mu)$ is the limiting probability distribution of a pure noise random matrix corresponding to the density $h(\alpha)$. Therefore, the pulled up value is the solution of

$$(5.15) \qquad \lambda - \alpha_1 = c\alpha_1 \int \frac{\mu}{\lambda - \mu} \, dF(\mu).$$

To finish the proof, we use Lemma 4 below, which shows that for any value of $c$ and density $h$, $\bar{m}(\lambda(\alpha)) = -1/\alpha$, and then insert this relation into the inverse equation (5.11). $\square$

LEMMA 4. *Let $\lambda(\alpha)$ denote the pulled up value corresponding to an original eigenvalue $\alpha$. For $\alpha$ large enough, regardless of the constant $c$ and of the underlying density $h(t)$,*

$$(5.16) \qquad \bar{m}(\lambda(\alpha)) = -\frac{1}{\alpha}.$$

PROOF. We rewrite (5.15) as follows:

$$\lambda(\alpha) = \alpha(1 - c) + \lambda \alpha c \int \frac{1}{\lambda - \mu} \, dF(\mu). \tag{5.17}$$

By extending the definition of the Stieltjes transform $m(z)$ of $F$, originally defined only for $z \in \mathbb{C}^+$, also to $z \in \mathbb{R}$ with $z > \text{support}(F)$, the last equation reads

$$\lambda(\alpha) = \alpha(1 - c) - \lambda \alpha c m(\lambda(\alpha)) \tag{5.18}$$

or stated otherwise

$$cm(\lambda(\alpha)) - \frac{1 - c}{\lambda(\alpha)} = -\frac{1}{\alpha}, \tag{5.19}$$

but according to (5.10), the left-hand side is simply $\bar{m}(\lambda(\alpha))$, also extended to the case of inputs $z \in \mathbb{R}$ with $z > \text{support}(F)$.  □

PROOF OF THEOREM 2.5. We consider the relation $\alpha(\lambda)$. That is, for each $\lambda > \text{support}(F)$ we look for the corresponding $\alpha$ such that (2.23) is satisfied. We show that there exists a unique solution $\alpha(\lambda)$, which is monotonic in $\lambda$, and as $\lambda \to \text{support}(F)$, satisfies that $\alpha(\lambda) \to \alpha^* < \infty$, but $d\alpha/d\lambda \to \infty$.

To this end, we rewrite (5.17) as follows:

$$\alpha = \frac{1}{(1 - c)/\lambda + c \int 1/(\lambda - \mu) \, dF(\mu)}. \tag{5.20}$$

Since for $\lambda > \text{support}(F)$ the functions $1/\lambda$ and $\int 1/(\lambda - \mu) \, dF(\mu)$ are both strictly monotonically decreasing in $\lambda$, it follows that

$$\frac{d\alpha}{d\lambda} > 0 \qquad \text{for } \lambda > \text{support}(F). \tag{5.21}$$

We now analyze the behavior of both $\alpha(\lambda)$ and its derivative as $\lambda$ approaches the support of $F$. According to [9, 34] near the boundary $b = \text{support}(F)$, the density $F$ exhibits a behavior closely resembling $\sqrt{|b - x|}$. Therefore,

$$\lim_{\lambda \to \text{support}(F)} \int \frac{1}{\lambda - \mu} \, dF(\mu) = \text{Const}, \tag{5.22}$$

whereas

$$\lim_{\lambda \to \text{support}(F)} \frac{d}{d\lambda} \int \frac{1}{\lambda - \mu} \, dF(\mu) = \infty. \tag{5.23}$$

This proves Theorem 2.5.  □

PROOF OF COROLLARY 3. We start our analysis from the equation

$$\lambda(\alpha) = \alpha(1 - c) + \alpha^2 c \int \frac{h(t)}{\alpha - t} \, dt. \tag{5.24}$$

First we make a change of variables $t = \mu_1 + s$, where

$$\mu_1 = \int t h(t) \, dt \tag{5.25}$$

and denote $h_0(s) = h(t)$, with the subscript zero denoting the fact that this density has zero mean. Then,

$$\lambda(\alpha) = \alpha(1-c) + \frac{\alpha^2 c}{\alpha - \mu_1} \int \frac{h_0(s)}{1 - s/(\alpha - \mu_1)} \, ds. \tag{5.26}$$

As $c \to \infty$, both $\alpha \to \infty$ and $\lambda \to \infty$. Specifically, $\alpha - \mu_1 > \text{support}(F)$. In this case we can expand

$$\frac{1}{1 - s/(\alpha - \mu_1)} = \sum_{k=0}^{\infty} \left( \frac{s}{\alpha - \mu_1} \right)^k,$$

where the sum is convergent for $|s| < \text{support}(F)$. Inserting this expansion into (5.24) gives

$$\lambda(\alpha) = \alpha(1-c) + \frac{\alpha^2 c}{\alpha - \mu_1} \left[ 1 + \int \frac{s^2 h_0(s)}{(\alpha - \mu_1)^2} \, ds + O\left( \frac{1}{\alpha - \mu_1} \right)^3 \right]. \tag{5.27}$$

Taking only the first two terms in the expansion gives the approximate solution

$$\alpha^* = \mu_1 (1 + \sqrt{c}).$$

To obtain the correction due to the variance of the distribution, we expand

$$\alpha^* = \mu \left( \beta_1 \sqrt{c} + \beta_0 + \frac{\beta_{-1}}{\sqrt{c}} + \cdots \right) \tag{5.28}$$

and insert the expansion into (5.27). This gives

$$\alpha^* = \mu \left( \sqrt{c} \sqrt{1 + \frac{\mu_2^2}{\mu_1^2}} + 1 + 4 \frac{\mu_2^2 / \mu_1^2}{1 + 2\mu_2^2 / \mu_1^2} + O\left( \frac{1}{\sqrt{c}} \right) \right) \tag{5.29}$$

for the location of the phase transition, and

$$\lambda(\alpha^*) = \mu \left[ c + 2\sqrt{c} \sqrt{1 + \frac{\mu_2^2}{\mu_1^2}} + O(1) \right]. \tag{5.30}$$

$\square$

5.4. *The phase transition phenomenon for finite $p, n$.* We conclude by presenting two examples of the phase transition phenomenon for the standard spiked covariance model with finite $p, n$. In Figure 1 we present an example of this phenomenon with $n = 50, p = 200, \kappa = 2.8$. For small noise level $\sigma$, the largest eigenvalue is roughly $\kappa^2 = 7.87$, and $\langle \mathbf{v}_{\text{PCA}}, \mathbf{e}_1 \rangle \approx 1$. As the noise level $\sigma$ increases the dot product decreases smoothly. However, at a noise level of roughly $\sigma = 1.85$ the largest eigenvalue and the second largest eigenvalue "cross" each other, leading to a *sharp* decrease in the overlap between the first eigenvector of PCA and the correct direction $\mathbf{e}_1$. The reason is that at this level of noise the spectral norm of the noise matrix overcomes the eigenvalue corresponding to the original signal. Thus, after the crossing the corresponding leading eigenvector is due to noise and points to a random direction.

As a second example, we present a possible behavior of $R = |\langle \mathbf{v}_{\mathrm{PCA}}, \mathbf{e}_1 \rangle|$ as a function of the number of samples $n$. In Figure 2 we present a simulation result with $p = 600$, $\sigma = 1$, $\|\mathbf{v}\| = 2$ and a variable number of samples $n$. When $n$ is very small the signal information is too weak, and as expected the overlap $R$ is very small. Then, as $n$ increases the largest eigenvalue of the noisy covariance matrix corresponds to the correct one and there is an increase in $R$. However, in this example at $n = 77$ there is a short crossover between the dominant eigenvalue and one due to noise (e.g., $\lambda_2 \approx \lambda_1$), leading first to a sharp decrease in $R$ and then to a sharp increase in $R$ back to around 0.5. The formulas derived in this paper can be used to estimate the probability of such a loss of locking to occur. This analysis also shows that care should be taken when applying bootstrap procedures for the eigenvectors, since in the case of weak signals in certain subsets of samples the resulting leading eigenvector might be due to noise.

## APPENDIX: PROBABILISTIC BOUNDS ON THE NORM OF WISHART MATRICES

Let $\Gamma$ denote an $n \times p$ matrix whose entries are all i.i.d. $N(0,1)$ random variables. Consider the $p \times p$ scaled Wishart matrix $W = (1/n)\Gamma^T\Gamma$ and the matrix $E = W - I_p$. In this Appendix we provide a probabilistic bound on the spectral norm of the matrices $W$ and $W - I_p$.

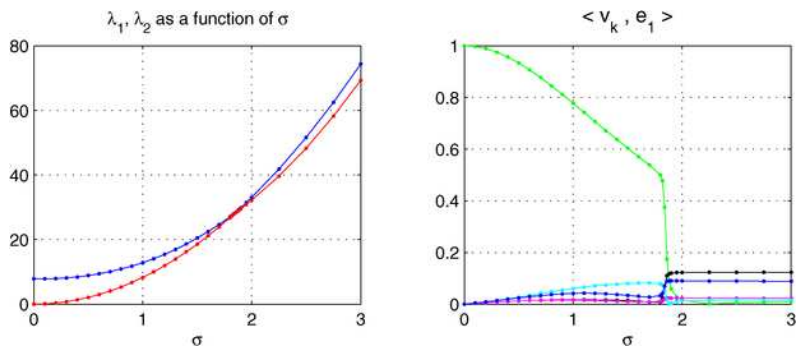

Fig. 1.  *Loss of tracking of the correct direction as a function of noise level.*

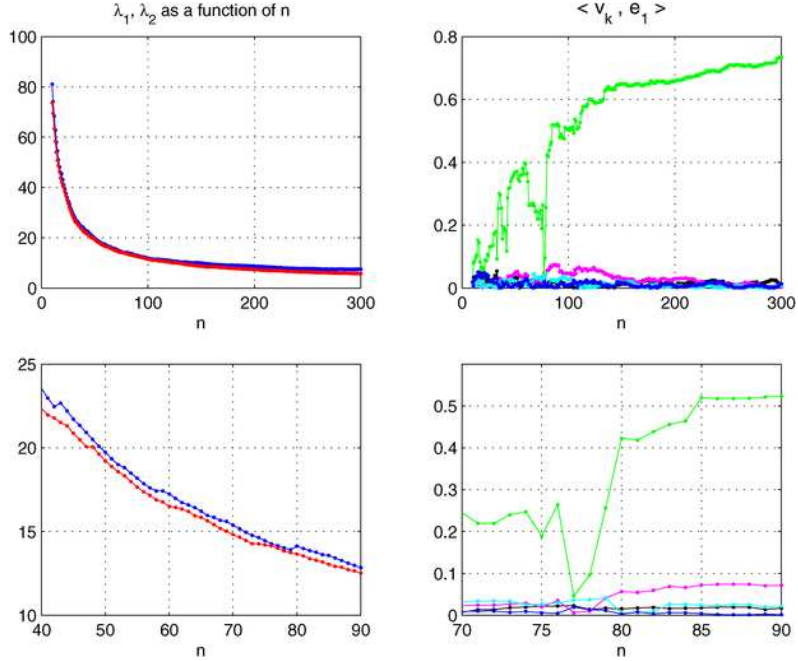

Fig. 2. *Loss of tracking of the correct direction as a function of sample size n. Bottom graphs are zoomed-in versions of the top ones.*

As shown in [18], in the limit $n, p \to \infty$, $p/n = c < 1$, the largest eigenvalue $\lambda_W$ of the matrix $W$ is distributed according to

$$\lambda_W = \frac{1}{n}\left[(\sqrt{n} + \sqrt{p-1})^2 + (\sqrt{n} + \sqrt{p-1})\left(\frac{1}{\sqrt{n}} + \frac{1}{\sqrt{p-1}}\right)^{1/3} W_1\right],$$

where $W_1$ follows a Tracy–Widom distribution of order 1. Therefore, we seek a bound on the norm of $E$ of the form $\|E\| \leq (1 + \sqrt{c})^2 + c - 1 + \text{const}$, that will hold with probability $1 - \varepsilon$. To prove such a bound (Lemma 2) we use the following proposition by Davidson and Szarek [7], which holds for finite $p, n$, with $n \leq p$.

PROPOSITION (Theorem II.13 in [7]). *Let $\Gamma$ be an $n \times p$ matrix with $n \leq p$ whose entries are all i.i.d. $N(0, 1)$ Gaussian variables. Let $s_1(\Gamma)$ be the largest singular value of $\Gamma$; then*

$$(A.1) \qquad \Pr\{s_1(\Gamma) > \sqrt{n} + \sqrt{p} + \sqrt{p}t\} < \exp(-pt^2/2).$$

PROOF OF LEMMA 2. The largest eigenvalue of $W = (1/n)\Gamma^T\Gamma$ is given by $\lambda_W = [s_1(\Gamma)]^2/n$. Therefore,

$$\Pr\{\lambda_W > (1 + \sqrt{p/n})^2 + \alpha p/n\} = \Pr\{s_1(\Gamma) > \sqrt{p}\sqrt{(1 + \sqrt{n/p})^2 + \alpha}\}.$$

We write

$$\sqrt{p}\sqrt{(1+\sqrt{n/p})^2+\alpha} = \sqrt{n}+\sqrt{p}+\sqrt{p}t$$

with

$$t = \sqrt{(1+\sqrt{n/p})^2+\alpha} - (1+\sqrt{n/p})$$

and use (A.1) to obtain that

$$\Pr\{\lambda_W > (1+\sqrt{p/n})^2 + \alpha p/n\}$$

(A.2)
$$\leq \exp\left[-\frac{p}{2}(\sqrt{\alpha+(1+\sqrt{n/p})^2}-(1+\sqrt{n/p}))^2\right]$$

$$\leq \exp\left[-\frac{p}{2}\left(\frac{\alpha}{\sqrt{\alpha+(1+\sqrt{n/p})^2}+1+\sqrt{n/p}}\right)^2\right].$$

We specifically consider $\alpha = 1$. Then, since $n \leq p$,

(A.3)
$$\frac{1}{\sqrt{1+(1+\sqrt{n/p})^2}+1+\sqrt{n/p}} \geq \frac{1}{\sqrt{5}+2}$$

and

(A.4)     $$\Pr\{\lambda_W > (1+\sqrt{p/n})^2 + p/n\} \leq \exp\left\{-\frac{p}{2(\sqrt{5}+2)^2}\right\} = \varepsilon.$$

Similarly, for $n/p \leq 1$, with probability at least $1 - \varepsilon$,

(A.5)     $$\|E\| = \|W - I\| \leq (1+\sqrt{p/n})^2 + \frac{p}{n} - 1 \leq 4\frac{p}{n}.$$     □

REMARK.   For $n < p$ the matrix $\Gamma^T\Gamma$ always has $p - n$ eigenvalues equal to 0, therefore $E$ has $p - n$ eigenvalues equal to $-1$. Thus, to bound $\|E\|$ we need only bounds on the largest positive eigenvalue as analyzed above.

**Acknowledgments.**   It is a pleasure to thank Iain Johnstone, Andrew Barron, Ilse Ipsen, Ann Lee and Patrick Perry for interesting discussions and useful suggestions. The author also acknowledges the anonymous referees for valuable criticism which greatly improved the exposition of the paper.

DEPARTMENT OF COMPUTER SCIENCE
AND APPLIED MATHEMATICS
WEIZMANN INSTITUTE OF SCIENCE
REHOVOT 76100
ISRAEL
E-MAIL: boaz.nadler@weizmann.ac.il